\documentclass[12pt,a4paper]{article}
\usepackage{amssymb,geometry,stmaryrd,color}
\usepackage{amssymb,amsfonts,amsmath,amsthm}
\usepackage{mathptmx}
\usepackage[doublespacing]{setspace}
\usepackage[english]{babel}
\usepackage[round,sort,comma,numbers]{natbib}
\usepackage{authblk}

\newtheorem{theorem}{\large Theorem}[section]

\newtheorem{definition}{\textsf{Definition}}
\newtheorem{remark}{\textsf{Remark}}

\usepackage{authblk}
\title{Fractional  non-homogeneous Poisson and P\'olya-Aeppli processes of order $ k $ and beyond}
\author[1]{Tetyana  Kadankova}
\author[2]{Nikolai Leonenko}
\author[3]{Enrico Scalas}
\affil[1]{Vrije Universiteit Brussel, Department of Mathematics, Belgium, Tetyana.Kadankova@vub.be}
\affil[2]{Cardiff University, School of Mathematics, UK, LeonenkoN@cardiff.ac.uk}
\affil[3]{University of Sussex, Department of Mathematics, UK, E.Scalas@sussex.ac.uk}
\date{}
\input{tcilatex}

\begin{document}

\maketitle

\noindent \textbf{Key words:} non-homogeneous fractional Poisson process of
order $k; $ non-homogeneous fractional P\'olya-Aeppli process of order $k; $
long range dependence; Caputo fractional derivative; $\alpha $-stable L\'evy
subordinators; fractional integro-differential difference equations.\newline
\noindent \textbf{AMS classification:} 60G55, 26A33, 60G05, 60G51

\noindent

\begin{center}
\textbf{Abstract}
\end{center}

\noindent We introduce two classes of point processes: a fractional
non-homogeneous Poisson process of order $k$ and a fractional
non-homogeneous P\'{o}lya-Aeppli process of order $k.$ 
We characterize these processes by deriving their non-local governing
equations. We further study the covariance structure of the processes and
investigate the long-range dependence property.

\section{Introduction}

Fractional Poisson processes (FPP) enjoy the property of non-stationarity
and long range dependence, which makes them an attractive modeling tool.
These processes are widely used in statistics, finance, meteorology, physics
and network science, see for instance \citep{Scalas} p. 332, and \citep{kuma}%
. \newline
Fractional Poisson processes were introduced as renewal processes in %
\citep{Mainardi1}. The authors generalized the characterization of the
Poisson process as the counting process for epochs defined as sum of
independent non-negative exponential random variables, and, instead of the
exponential, the authors used a Mittag-Leffler distribution. The theory of
FPP was further developed by Beghin and Orsingher \citep{Beghim1, Beg} and
by Meerschaert \emph{et al.} \citep{Meer3}. \newline
In particular, Meerschaert \emph{et al.} (2011) defined FPP by means of a
time-change for the Poisson process $N(t)$, where the time variable $t$ is
replaced by the inverse $\alpha $-stable subordinator $Y_{\alpha }(t)$.
Remarkably, they could prove the equality in distribution between $%
N(Y_{\alpha }(t))$ and the counting process defined by \citep{Mainardi1}. 
\newline
Leonenko \emph{et al.} \citep{Leo} used the same time-change technique to
introduce a non-homogeneous fractional Poisson process (NFPP) by replacing
the time variable in the FPP with an appropriate function of time. 
\newline
In a recent paper Gupta \emph{et al.} \citep{gupta1,gupta2} generalize the
results available on fractional Poisson processes using the $z$-transform
technique. \newline
Kostadinova and Minkova \citep{kost} introduced a Poisson process of order $%
k $ with insurance modelling in mind. This process models the claim arrival
in groups of size $k$, where the number of arrivals in a group is uniformly
distributed over $k$ points. \newline
The P\'{o}lya-Aeppli process of order $k$ was studied in \citep{Chuko2015}
and later by \citep{kost2019}. In this process, the uniform distribution on
the integers $1,\ldots ,k$ is replaced by the truncated geometric
distribution of parameter $\varrho $.\newline
To deal with dependent inter-arrival times, a generalization of Poisson
processes of order $k$ was proposed by Sengar \emph{et al.} \citep{Mahe1}.
These authors extended the Poisson process of order $k$ by means of time
change with a general L\'{e}vy subordinator as well as an inverse L\'{e}vy
subordinator. 
\newline
Here, we combine the compound Poisson processes of order $k$ and fractional
Poisson processes, namely we study a fractional non-homogeneous Poisson
process of order $k$ and a fractional non-homogeneous P\'{o}lya-Aeppli
process of order $k$ (see the definitions below). First, we generalize the
results of \citep{kost} by considering a non-homogeneous Poisson process of
order $k$. Then, we generalize the results of Sengar \emph{et al.} %
\citep{Mahe1} by introducing the time non-homogeneity in the fractional
Poisson process of order $k.$ Finally, we study a non-homogeneous fractional
P\'{o}lya-Aeppli process of order $k.$ \newline
This paper is organized as follows. Section 2 collects some known results
from the theory of subordinators and provides the definition of the compound
distributions of order $k.$ 
In Section 3, we consider a non-homogeneous fractional Poisson process of
order $k$. We obtain the governing equations and calculate the moments and
the covariance function of the process. Section 4 is devoted to a
non-homogeneous fractional Polya-Aeppli process of order $k.$ We derive the
non-local governing equations for the marginal distributions of these
processes, using non-local operators known as Caputo derivatives. The
moments and the covariance structure of the processes are derived, as well.
We conclude in Sections 5 by presenting some simulations and discussing
possible applications and generalizations. 


\section{Preliminaries}

This section presents known results in the theory of subordinators and
provides the definition of the compound distributions of order $k$.

\subsection{\textit{Compound distributions of of order $k$ }}

Consider a random variable that can be represented as a random sum $%
N=X_{1}+X_{2}+\ldots +X_{Y},$ where $\{X_{i}\}_{i=1}^{\infty }$ is a
sequence of independent identically distributed random variables (i.i.d.
r.v's), independent of a non-negative integer-valued random variable $Y$.
The probability distribution of $N$ is called compound distribution and the
distribution of $X_{1}$ is called compounding distribution. \newline
A well-known and widely used example is the compound Poisson distribution,
where $Y$ has a Poisson distribution. If $X_{i}\in \{1,2,\dots,k\}$, then
the random variable $N$ has a compound discrete distribution of order $k.$

\noindent Compound discrete distributions of order $k$ were studied by
Philippou \citep{phil1983} and \citep{Phil}. \newline
\noindent As mentioned previously, in this paper we deal with two types of
compounding distributions: the discrete uniform distribution and the
truncated geometric distribution. They respectively induce the Poisson
distribution of order $k$ and the P\'{o}lya-Aeppli distribution of order $k$
as will be shown in the following. We say that the random variable $X$ is
uniformly distributed over $k$ points if its probability mass function (pmf)
is of the form 
\begin{equation}
\mathbb{P}[X=m]=\frac{1}{k},\quad m=1,\dots ,k.  \label{unif}
\end{equation}%
Its probability generating function (pgf) 
\begin{equation*}
G_{X}(u)=\mathbb{E}\left[ u^{X}\right] =\frac{1}{k}(u+u^{2}+...+u^{k})=\frac{%
u}{k}\cdot \frac{1-u^{k}}{1-u}, \quad u\in (0,1).
\end{equation*}
\noindent The random variable $X$ has a truncated geometric distribution
with parameter $\varrho $ and with success probability $1-\varrho $ if 
\begin{equation}
\mathbb{P}[X=m]=\frac{1-\varrho }{1-\varrho ^{k}}\varrho ^{m-1},\qquad
m=1,2,\dots ,k,\qquad \varrho \in \lbrack 0,1).  \label{TG}
\end{equation}%
Consequently, the pgf of $X$ is given by 
\begin{equation}
G_{X}(u)=\mathbb{E}\left[ u^{X}\right] =\frac{(1-\varrho )u}{1-\varrho ^{k}}:%
\frac{1-\varrho ^{k}u^{k}}{1-\varrho u}, \quad u\in (0,1).  \label{pgf}
\end{equation}%
Note, that for $k\rightarrow \infty,$ the truncated geometric distribution
asymptotically coincides with the geometric distribution with parameter $%
1-\varrho.$\newline
We can now define the Poisson distribution of order $k$.

\begin{definition}[Poisson distribution of order $k$]
The random variable $N$ has Poisson distribution of order $k$ \ with
parameter $\Lambda $ if $N=X_{1}+X_{2}+\dots +X_{Y},$ where: \newline
(1) $\{X_{i}\}_{i\geq 1}$ are the i.i.d. r.v's with the uniform
distribution; (2) $Y$ has Poisson distribution with parameter $\Lambda >0; $
(3) $Y$ and $\{X_{i}\}_{i\geq 1}$ are independent.
\end{definition}

\noindent Note that%
\begin{equation*}
\mathbb{P}[N=m]=e^{-\Lambda k}\sum_{(n_{1},\dots n_{k})\in \Omega (k,m)}%
\frac{\Lambda ^{n_{1}+...+n_{k}}}{n_{1}!\cdot ...\cdot n_{k}!}=e^{-\Lambda
k}\sum\limits_{\Omega (k,m)}\frac{\Lambda ^{z_{k}}}{\Pi _{k}!},
\end{equation*}%
where $z_{k}=n_{1}+n_{2}+\dots +n_{k},\quad \Pi _{k}!=n_{1}!\cdot
n_{2}!\cdot \dots \cdot n_{k}!,$ and 
\begin{equation}
\Omega (k,m)=\{(n_{1},\dots n_{k}):n_{1}+2n_{2}+\dots kn_{k}=m\}.
\label{OMEGA}
\end{equation}%
The pgf of the Poisson distribution of order $k$ 
\begin{equation}
G_{N}(u)=\mathbb{E}\left[ u^{N}\right] =\exp \left\{ -\Lambda \left( k-\sum
\limits_{j=1}^{k}u^{j}\right) \right\} .  \label{MGF1}
\end{equation}%
\noindent Note that $N\overset{d}{=}\sum \limits_{j=1}^{k}jY_{j},$ where $%
Y_{j},j=1,\dots ,k$ are independent copies of Poisson random variable $Y$
with parameter $\Lambda $, and \textquotedblleft $\overset{d}{=}$%
\textquotedblright\ stands for equality in distributions. \noindent We now
introduce the P\'{o}lya-Aeppli distribution of order $k$ as a compound
Poisson distribution (see \citep{Min}).

\begin{definition}[P\'{o}lya-Aeppli distribution of order $k$]
The random variable $N$ has P\'{o}lya-Aeppli distribution of order $k$ with
parameter $1-\varrho $ if $N=X_{1}+X_{2}+\dots +X_{Y},$ where: (i) $%
\{X_{i}\}_{i\geq 1}$ are the i.i.d. r.v's with the truncated geometric
distribution with parameter $1-\varrho $, given by (\ref{TG}); (ii) $Y$ has
Poisson distribution with parameter $\Lambda;$ (iii) $Y$ and $%
\{X_{i}\}_{i\geq 1}$ are independent.
\end{definition}

\noindent Note that the probability generating function of $N$ is $G_{N}(u)=%
\mathrm{e}^{-\Lambda (1-P_{X_{1}}(u))},$ where $P_{X_{1}}$ is given by (\ref%
{pgf}).

The probability mass function of P\'{o}lya-Aeppli distribution of order $k$
is defined by (see \citep{Min}, Theorem 3.1):

\begin{equation}
\mathbb{P}[N=m]=q_{m}(\Lambda ),m=0,1,2,...,  \label{PAPMF}
\end{equation}

where 
\begin{equation*}
q_{0}(\Lambda )=e^{-\Lambda }
\end{equation*}%
\begin{equation*}
q_{m}(\Lambda )=e^{-\Lambda }\sum_{j=1}^{m}{\binom{m-1}{j-1}}\frac{Q^{J}}{j!}%
\rho ^{m-j},m=1,2,...,k
\end{equation*}%
\begin{equation*}
q_{m}(\Lambda )=e^{-\Lambda }\left[\sum_{j=1}^{m}{\binom{m-1}{j-1}}\frac{Q^{J}}{j!%
}\rho ^{m-j}-\sum \limits_{n=1}^{l}(-1)^{n-1}\frac{(Q\rho ^{k})^{n}}{n!}%
\times \right.
\end{equation*}

\begin{equation*}
\left. \times \sum \limits_{j=0}^{m-n(k+1)}{\binom{m-n(k+1)+n-1}{j+n-1}}\frac{Q^{J}%
}{j!}\rho ^{m-j-n(k+1)}\right],
\end{equation*}

and 
\begin{equation*}
Q=\frac{\Lambda (1-\rho )}{1-\rho ^{k}},m=l(k+1)+r,r=0,1,..,k,l=1,2,....
\end{equation*}

\subsection{Inverse $\protect\alpha$-stable subordinator}

Let $\mathcal{L}_{\alpha }=\{\mathcal{L}_{\alpha }(t); t\geq 0\}$ be a $%
\alpha $-stable L\'{e}vy subordinator, that is L\'{e}vy process with Laplace
transform: 
\begin{equation*}
\mathbb{E} \left[ \mathrm{e}^{-s\mathcal{L}_{\alpha }(t)} \right]%
=e^{-ts^{\alpha }},\quad 0<\alpha <1,\quad s\geq 0.
\end{equation*}%
Then the inverse $\alpha $-stable subordinator $\{Y_{\alpha }(t); \: t\geq
0\}$ (see e.g. \citep{Meer} p. 103) is defined as the first passage time of $%
\mathcal{L}_{\alpha }:$ 
\begin{equation}  \label{inverse}
Y_{\alpha} (t) =\inf\{u>0: \mathcal{L}_{\alpha}(u) > t \}, \quad t \ge 0.
\end{equation}
We will use the following properties of the inverse $\alpha $-stable
subordinator:

\begin{itemize}
\item[(i)] The density of $Y_{\alpha }(t)$ is of the form (see \citep{Meer}
p.113): 
\begin{equation}  \label{densityinv}
h_{\alpha }(t,x)=\frac{d}{dx}\mathbb{P}[Y_{\alpha }(t)\leq x]=\frac{t}{%
\alpha }x^{-1-\frac{1}{\alpha }}g_{\alpha }(tx^{-\frac{1}{\alpha }}),\quad
x>0,\quad t>0,
\end{equation}%
where 
\begin{equation*}
g_{\alpha }(x)=\frac{1}{\pi }\sum\limits_{k=1}^{\infty }(-1)^{k+1}\frac{%
\Gamma (\alpha k+1)}{k!}\frac{1}{x^{\alpha k+1}}\sin (\pi k\alpha )
\end{equation*}%
is the density of $\mathcal{L}_{\alpha }(1)$ (see e.g. \citep{kataria}).

\item[(ii)] The Laplace transform 
\begin{equation}
\tilde{h}_{\alpha }(s,x)=\int\limits_{0}^{\infty }\mathrm{e}^{-st}h_{\alpha
}(t,x)dt=s^{\alpha -1}e^{-xs^{\alpha }},\quad s\geq 0.  \label{LTIS}
\end{equation}

\item[(iii)] The moments of the inverse $\alpha $-stable subordinator are as
follows: 
\begin{equation}
\mathbb{E}[Y_{\alpha }^{\nu }(t)]=\frac{\Gamma (\nu +1)}{\Gamma (\alpha \nu
+1)}t^{\alpha \nu },\nu >0,\quad \mathrm{Var}[Y_{\alpha }(t)]=t^{2\alpha }%
\left[ \frac{2}{\Gamma (2\alpha +1)}-\frac{1}{(\Gamma (\alpha +1))^{2}}%
\right] .  \label{var}
\end{equation}%
(see e.g. \citep{kataria} p.1640).

\item[(iv)] The covariance function (see \citep{Leo1,Leo}) is 
\begin{equation}
\mathrm{Cov}[Y_{\alpha }(t),Y_{\alpha }(s)]=\frac{1}{\Gamma (1+\alpha
)\Gamma (\alpha )}\int\limits_{0}^{\min (t,s)}((t-\tau )^{\alpha }+(s-\tau
)^{\alpha })\tau ^{\alpha -1}d\tau -\frac{(st)^{\alpha }}{\Gamma
^{2}(1+\alpha )}.  \label{covinv}
\end{equation}
\end{itemize}

\section{ Poisson processes of order $k$}

The Poisson process of order $k$ was introduced in \citep{kost}, see also %
\citep{Mahe1}.

\begin{definition}
The Poisson process of order $k$ (PPk) $N=\{N(t); \: t\geq 0\}$ is defined
as a compound Poisson process with the compounding discrete uniform
distribution: 
\begin{equation}
N(t)=X_{1}+\dots +X_{N_{1}(t)},  \label{PPk}
\end{equation}%
where (1) $X_{i}$ are independent copies of a discrete uniform random
variable distributed over $k$ points given by (\ref{unif}); (2) $%
N_{1}=\{N_{1}(t); \: t\geq 0\}$ is the Poisson process with parameter $%
k\lambda;$ (3) $N_{1}$ and $\{X_{i}\}_{i\geq 1}$ are independent.
\end{definition}

\noindent The following Kolmogorov forward equations are valid for $p_{m}(t)=%
\mathbb{P}[N(t)=m]:$ 
\begin{align}
& \frac{d}{dt}p_{0}(t)=-k\lambda p_{0}(t)  \label{govern_PPk} \\
& \frac{d}{dt}p_{m}(t)=-k\lambda p_{m}(t)+\lambda \sum\limits_{j=1}^{m\wedge
k}p_{m-j}(t),\quad m=1,2,\dots
\end{align}%
with the initial condition $p_{0}(0)=1,:p_{m}(0)=0,m\geq 1,$ and $m\wedge
k=\min (m,k).$ The pgf is of the form: 
\begin{equation*}
G_{N(t)}(u)=\mathbb{E}\left[ u^{N(t)}\right] =\exp \{\lambda
t(u+...+u^{k}-k)\},
\end{equation*}%
and the first two moments are given by 
\begin{equation}
\mathbb{E}[N(t)]=\frac{k(k+1)}{2}\lambda t,\quad \mathrm{Cov}[N(t),N(s)]=%
\frac{k(k+1)(2k+1)}{6}\lambda \min (s,t).  \label{VARN}
\end{equation}

\subsection{\textit{\ Fractional Poisson process of order $k $}}

In this sub-section we shall derive governing equations for a fractional
Poisson process of order $k$ and we shall investigate its long-range
dependence properties. It is worth noting that Sengar et al. \citep{Mahe1}
studied the Poisson process of order $k$ time-changed by a general L\'{e}vy
subordinator and its inverse. However, among their examples, they did not
explicitly consider the governing equations for the inverse $\alpha $-stable
subordinator (this particular process is studied in \citep{gupta2}). That is
why we specify some formulae of \citep{Mahe1} that will be used in the
following sub-sections. In particular, below, we use equation (\ref{var}) to
derive the marginal distributions of the fractional Poisson process of order 
$k$. 

\begin{definition}
(Fractional Poisson process of order $k$). The process $N_{\alpha }(t)$ is
called fractional Poisson process of order $k$ (FPPk) if 
\begin{equation}
N_{\alpha }(t)=N(Y_{\alpha }(t)),\quad \quad 0<\alpha <1,  \label{FPPk}
\end{equation}%
where (1) $Y_{\alpha} (t)$ is the inverse $\alpha$-stable subordinator,
given by (\ref{inverse}); (2) $N $ is the Poisson process of order $k,$
given by (\ref{PPk}); (3) $Y_{\alpha} (t)$ and $N$ are independent. 
\end{definition}

\noindent The marginal distributions of the FPPk process is given by 
\begin{equation*}
p_{m}^{\alpha }(t)=\mathbb{P}[N_{\alpha }(t)=m]=\sum\limits_{\Omega (k,m)}%
\frac{\lambda ^{z_{k}}}{\Pi _{k}!}\sum\limits_{n=0}^{\infty }\frac{%
(-k\lambda )^{n}}{n!}\mathbb{E}\left[ (Y_{\alpha }(t))^{z_{k}+n}\right] =
\end{equation*}%
\begin{equation*}
=\sum\limits_{\Omega (k,m)}\frac{\lambda ^{z_{k}}}{\Pi _{k}!}%
\sum\limits_{n=0}^{\infty }\frac{(-k\lambda )^{n}}{n!}\frac{\Gamma
(z_{k}+n+1)}{\Gamma (\alpha (z_{k}+n)+1)}t^{\alpha (z_{k}+n)},m=0,1,..
\end{equation*}%
where $z_{k}=n_{1}+n_{2}+\dots +n_{k},\quad \Pi _{k}!=n_{1}!n_{2}!\dots
n_{k}!$, and $\ \Omega (k,m)$ is defined in (\ref{OMEGA}).

Also%
\begin{equation*}
\mathbb{E}[N_{\alpha }(t)]=k\lambda \frac{(k+1)}{2}\mathbb{E}[Y_{\alpha
}(t)],
\end{equation*}%
\begin{equation*}
\mathrm{Var}[N_{\alpha }(t)]=k\lambda \frac{(k+1)(2k+1)}{6}\mathbb{E}%
[Y_{\alpha }(t)]+\left( k\lambda \frac{(k+1)}{2}\right) ^{2}\mathrm{Var}%
(Y_{\alpha }(t)),
\end{equation*}%
\begin{equation*}
\mathrm{Cov}[N_{\alpha }(t),N_{\alpha }(s)]=\frac{k(k+1)(2k+1)\lambda (\min
(t,s))^{\alpha }}{6\Gamma (1+\alpha )}+\left( \frac{k\lambda (k+1)}{2}%
\right) ^{2}\mathrm{Cov}(Y_{\alpha }(s),Y_{\alpha }(t)),
\end{equation*}%
where the variance and covariance of the process $Y_{\alpha }(t)$ are given
by (\ref{var}) and (\ref{covinv}). 

\subsection*{\textbf{Correlation structure and long-range dependence}}

\noindent There exist many definitions of the long-range dependence
property. Here, we shall use the definition given in \citep{Biard1}.

\begin{definition}
The process $\{X(t); \: t\geq 0\}$ has a long-range dependence property
(LRD) if for fixed $s$ and some $c(s)$ and $\alpha \in
(0,1):\lim\limits_{t\rightarrow \infty }[\mathrm{Corr}(X(s),X(t))/t^{-\alpha
}]=c(s),$ where $\mathrm{Corr}$ is the correlation function of the process $%
X $.
\end{definition}

\noindent We now investigate the asymptotic behavior of the correlation
function of the FPPk process defined by (\ref{FPPk}).

\begin{theorem}
The process $N_{\alpha}(t)$ has the LRD property.
\end{theorem}

\noindent \textit{\ Proof.} 
Using the result of \citep{Leo1} we have that for a fixed $s>0$ 
\begin{equation*}
\mathrm{Corr}[N_{\alpha }(t),N_{\alpha }(s)]\sim t^{-\alpha }C(\alpha
,s)\qquad t\rightarrow \infty ,
\end{equation*}%
where $C(\alpha ,s)=\left( \frac{1}{\Gamma (2\alpha )}-\frac{1}{\alpha
(\Gamma (\alpha ))^{2}}\right) ^{-1}\left[ \frac{\alpha \mathrm{Var}[N(1)]}{%
\Gamma (1+\alpha )(\mathbb{E}[N(1)])^{2}}+\frac{\alpha s^{\alpha }}{\Gamma
(1+2\alpha )}\right] , $ and $\mathbb{E}[N(1)]$ and $\mathrm{Var}[N(1)]$ are
given by (\ref{VARN}).\newline
\noindent \textbf{\ Governing equations}\newline
\noindent In the sequel we will employ the fractional Caputo (or
Caputo-Djrbashian) derivative which is defined as follows (see e.g. %
\citep{Meer} p. 30) 
\begin{equation}  \label{caputo}
D_{t}^{\alpha }f(t)=%
\begin{cases}
& \frac{1}{\Gamma (\alpha )}\int\limits_{0}^{t}\frac{df(u)}{du}\frac{du}{%
(t-u)^{\alpha }},\qquad 0<\alpha <1, \\ 
& \frac{df(u)}{du},\qquad \qquad \qquad \qquad \alpha =1.%
\end{cases}%
\end{equation}%
%
%
%
%
%
%
%
%
%
%
%
%
%
%
%
%
%

\begin{theorem}
The governing fractional difference-differential equations for \newline
$p_{m}^{\alpha }(t),t\geq 0$ are given by 
\begin{align}
& D_{t}^{\alpha }p_{0}^{\alpha }(t)=-k\lambda p_{0}^{\alpha }(t) \\
& D_{t}^{\alpha }p_{m}^{\alpha }(t)=-k\lambda p_{m}^{\alpha }(t)+\lambda
\sum\limits_{j=1}^{m\wedge k}p_{m-j}^{\alpha }(t),\quad m=1,2,\dots
\end{align}%
with the initial condition 
\begin{equation*}
p_{m}^{\alpha }(0)=\delta _{m,0}=%
\begin{cases}
1,\quad m=0 \\ 
0,\quad m\geq 1.%
\end{cases}%
\end{equation*}
\end{theorem}

\noindent Note, that by setting $\alpha =1,$ we get the governing equations
of the Poisson process of order $k$ given in equation (\ref{govern_PPk}). 
\newline
\noindent \textit{\ Proof.} Note that 
\begin{equation}
D_{t}^{\alpha }h_{\alpha }(t,u)=-\frac{\partial }{\partial u}h_{\alpha }(t,u)
\label{governinvsub}
\end{equation}%
and remember that 
\begin{equation}
p_{n}^{\alpha }(t)=\int\limits_{0}^{\infty }p_{n}(u)h_{\alpha }(t,u)du\qquad
n=0,1,2...  \label{convolutionfppk}
\end{equation}%
We first consider the case $n\geq 1.$ By taking the fractional Caputo
derivative of both sides (\ref{convolutionfppk}) and using property (\ref%
{governinvsub}), we get 
\begin{align*}
D_{t}^{\alpha }p_{m}^{\alpha }(t)& =-\int\limits_{0}^{\infty }p_{m}(u)\frac{%
\partial }{\partial u}h_{\alpha }(t,u)du= \\
& =\int\limits_{0}^{\infty }[-k\lambda p_{m}(u)+\lambda
\sum\limits_{j=1}^{m\wedge k}p_{m-j}(u)]h_{\alpha }(t,u)du-p_{m}(u)h_{\alpha
}(t,u)|_{0}^{\infty }= \\
& =-k\lambda p_{m}^{\alpha }(t)+\lambda \sum\limits_{j=1}^{m\wedge
k}p_{m-j}^{\alpha }(t).
\end{align*}

%
%
%
%
%
%
\noindent For $n=0$ we have%
\begin{equation*}
D_{t}^{\alpha }p_{0}^{\alpha }(t)=-\int\limits_{0}^{\infty }p_{0}(u)\frac{%
\partial }{\partial u}h_{\alpha }(t,u)du=\int\limits_{0}^{\infty }[-k\lambda
p_{0}(u)]h_{\alpha }(t,u)du=-k\lambda p_{0}^{\alpha }(t).
\end{equation*}

\begin{remark}
Note that Sengar et al. \citep{Mahe1} derived governing equations in which
the Caputo derivative is replaced by a more general non-local operator. We
present the proof of Theorem 3.3 for the sake of completeness.
\end{remark}

\subsection{ \textit{Non-homogeneous Fractional Poisson process of order $k $%
}}

We now generalize the fractional Poisson process of order $k$ by introducing
a deterministic, time dependent intensity or rate function $\lambda
(t):[0,\infty )\rightarrow \lbrack 0,\infty ),$ such that for every fixed $t
> 0,$ the cumulative rate function%
\begin{equation*}
\Lambda (t)=\int\limits_{0}^{t}\lambda (u)du<\infty
\end{equation*}%
Denote $\Lambda (s,t)=\int\limits_{s}^{t}\lambda (u)du=\Lambda (t)-\Lambda
(s), \:0\leq s < t.$ Let $N_{1}^{1}(t); \: t\geq 0$ be a homogeneous Poisson
process (HPP) of unit intensity, and $N_{1}^{1}(\Lambda (t)), \: t\geq 0,$
be a non-homogeneous Poisson process (NPP) with rate function $\lambda (t),$
then 
\begin{equation*}
N^{n}(t)=X_{1}+\ldots +X_{N_{1}^{1}(k\Lambda (t))}, \quad t\geq 0,
\end{equation*}%
is non-homogeneous Poisson process of order $k$ (NPPk), with rate function $%
\lambda (t), \: t\geq 0,$ where $\{X_{i}\}_{i\geq 1}$ are the i.i.d.r.v's
with the uniform distribution on $\{1,2, \dots, k\},$ independent of $%
N_{1}^{1}(\Lambda (t)).$ The mgf of the process $N^{n}$ is of the form: 
\begin{equation*}
G_{N^{n}(t)}(u)=\mathbb{E} \left[ u^{N^{n}(t)} \right] =\exp \{\Lambda
(t)(u+...+u^{k}-k)\}.
\end{equation*}
\noindent The process $N^{n}$ has the following distributions of its
increments:%
\begin{equation*}
p_{m}^{n}(t,u)=\mathbb{P}[N^{n}(t+u)-N^{n}(u)=m]=
\end{equation*}%
\begin{equation}
\mathbb{=}e^{-k\Lambda (u,t+u)}\sum_{\Omega (k,m)}\frac{[\Lambda
(u,u+t)]^{n_{1}+...+n_{k}}}{n_{1}!...n_{k}!}, \quad m=0,1, \dots
\label{nonhomk}
\end{equation}%
Incidentally, this model includes Weibull's rate function: $\Lambda
(t):=\Lambda (0,t)=\left( \frac{t}{b}\right) ^{c},\quad \lambda (t)=\frac{c}{%
b}\left( \frac{t}{b}\right) ^{c-1},\quad c\geq 0,b>0;$ Makeham's rate
function: $\Lambda (t)=\frac{c}{b}e^{bt}-\frac{c}{b}+\mu t,\quad \lambda
(t)=ce^{bt}+\mu ,\quad c>0,b>0,\mu \geq 0,$ and many others. \newline
\noindent We define a non-homogeneous fractional Poisson process of order $k$
(FNPPk) as 
\begin{equation}
N_{\alpha }^{\ast }(t)=N^{n}(Y_{\alpha }(t)),\quad t\geq 0,\quad 0<\alpha <1,
\label{FNPPk}
\end{equation}%
where $Y_{\alpha }(t)$ is the inverse $\alpha $-stable subordinator (\ref%
{inverse}), independent of the NPPk process $N^{n}.$ \newline
\noindent \textbf{Marginal distributions}\newline
\noindent Define the increment process: $I_{\alpha }(t,v)=N(\Lambda
(Y_{\alpha }(t)+v))-N(\Lambda (v)).$ 
Its marginal distributions can be written as follows: 
\begin{equation}
p_{m}^{\ast }(t,v)=\mathbb{P}[I_{\alpha }(t,v)=m]=\int\limits_{0}^{\infty
}p_{m}^{n}(u,v)h_{\alpha }(t,u)du,  \label{Main}
\end{equation}%
where $h_{\alpha }(t,u)$ is the density of the inverse $\alpha $-stable
subordinator (\ref{densityinv}) and $p_{x}^{n}(u,v)$ is given by (\ref%
{nonhomk}). Consequently the marginal distributions of $N_{\alpha }^{\ast
}(t)$ are given by 
\begin{equation*}
\mathbb{P}[N_{\alpha }^{\ast }(t)=m]=p_{m}^{\ast
}(t,0)=\int\limits_{0}^{\infty }p_{m}^{n}(u,0)h_{\alpha }(t,u)du.
\end{equation*}%
For the NFPP $N_{1}^{1}(\Lambda (Y_{\alpha }(t)); \: t\geq 0, $ of order $%
k=1,$ Leonenko \emph{et al.} \citep{Leo} derived the governing equations for
the marginal distributions $\mathbb{P}[I_{\alpha }^{1}(t,v)=m]$ of the
corresponding increment process $I_{\alpha }^{1}(t,v)=N_{1}^{1}(\Lambda
(Y_{\alpha }(t)+v))-N_{1}(\Lambda (v))$ of NFPP (of order $k=1$), where $%
N_{1}^{1}$ is the homogeneous Poisson process of intensity $1$. \noindent We
shall derive the governing equations for the marginal distributions $%
p_{x}^{\ast }(t,v)$ of FNPPk.

\begin{theorem}
The marginal distributions $p_{x}^{\ast }(t,v)$ satisfy the following
fractional differential-difference integral equations 
\begin{align}
& D_{t}^{\alpha }p_{0}^{\ast }(u,v)=-k\int\limits_{0}^{\infty }\lambda
(u+v)p_{0}^{n}(u,v)h_{\alpha }(t,u)du\quad \quad 0\leq v<u \\
& D_{t}^{\alpha }p_{m}^{\ast }(u,v)=\int\limits_{0}^{\infty }[-k\lambda
(u+v)p_{m}^{n}(u,v)+\lambda (u+v)\sum\limits_{j=1}^{m\wedge
k}p_{m-j}^{n}(u,v)]h_{\alpha }(t,u)du,\quad m=1,2,\dots  \notag
\end{align}%
with the initial condition: $p_{m}^{\ast }(0,v)=\delta _{m,0},$ where $%
p_{m}^{n}(u,v)$ is given by (\ref{nonhomk}).
\end{theorem}

\noindent \textit{\ Proof.} Note that the mgf of $p_{m}^{n}(u,v)$ is of the
form 
\begin{equation*}
\hat{p}_{s}^{n}(u,v)=\mathbb{E} \left[ s^{N^{n}(v+u)-N^{n}(v)} \right]=\exp
\{\Lambda (v,u+v)(s+\dots+s^{k}-k)\},
\end{equation*}
while the Laplace transform with respect to $t$ of $h_{\alpha }(t,u)$ is
given by (\ref{LTIS}). Taking both the mgf and the Laplace transform in (\ref%
{Main}) as above, we have 
\begin{equation}
\bar{p}_{s}^{\ast }(r,v)=\int\limits_{0}^{\infty }\hat{p}_{s}(u,v)\tilde{h}%
_{\alpha }(r,u)du=r^{\alpha -1}\int\limits_{0}^{\infty }\exp \{\Lambda
(v,u+v)(s+\dots+s^{k}-k)\}e^{-ur^{\alpha }}du.  \label{1}
\end{equation}
Note that for $U(u)=\exp \{\Lambda (v,u+v)(s+\dots+s^{k}-k)\},$ we have 
\begin{equation}
\frac{d}{du}U(u)=(s+s^{2}+...+s^{k}-k)\lambda (u+v)\exp \{\Lambda
(v,u+v)(s+s^{2}+...+s^{k}-k)\}.  \label{2}
\end{equation}
Thus, integrating (\ref{1}) by parts with $U$ as above, and $%
V=-e^{-ur^{\alpha }}/r^{\alpha }$, we get 
\begin{align}
& \bar{p}_{s}^{\ast }(r,v)=r^{\alpha -1}\left \{ \left[ -\frac{1}{r^{\alpha }%
}(\mathrm{e}^{-ur^{\alpha }}\mathrm{e}^{\Lambda
(v,u+v)(s+...s^{k}-k)}|_{0}^{\infty } \right] + \right.  \label{3} \\
+& \left. \frac{1}{r^{\alpha }}(s+s^{2}+...+s^{k}-k)\int\limits_{0}^{\infty
}k\lambda (v,u+v)\exp \{\Lambda (v,u+v)(s+s^{2}+...+s^{k}-k)\}\mathrm{e}%
^{-ur^{\alpha }}du \right \}=  \notag \\
=& \frac{1}{r^{\alpha }}\left[ r^{\alpha
-1}+(s+s^{2}+...+s^{k}-k)\int\limits_{0}^{\infty }\lambda (u+v)\exp
\{\Lambda (v,u+v)(s+s^{2}+...+s^{k}-k)\}r^{\alpha -1}\mathrm{e}^{-ur^{\alpha
}}du\right] .  \notag
\end{align}%
\noindent We shall use the following property of the Caputo derivative: 
\begin{equation*}
\mathcal{L}_{r}\{D_{t}^{\alpha }f\}(r)=r^{\alpha }\mathcal{L}%
\{f\}(r)-r^{\alpha -1}f(0^{+}),
\end{equation*}%
where $\mathcal{L}\{(f)\}(r)$ stands for the Laplace transform of function $%
f.$ Note that $p_{y}^{\ast }(0^{+},v)=1, $ since \: $Y_{\alpha }(0)=0$ a.s.
Hence, by (\ref{3}) 
\begin{align*}
& r^{\alpha }\bar{p}_{s}^{\ast }(r,v)-r^{\alpha -1}\bar{p}_{s}^{\ast }(0,v)=%
\mathcal{L}_{r}\{D_{t}^{\alpha }\bar{p}_{s}^{\ast }(r,v)\}(r)= \\
& =(s+s^{2}+...+s^{k}-k)\int\limits_{0}^{\infty }\lambda (u+v)\exp \{\Lambda
(v,u+v)(s+s^{2}+...+s^{k}-k)\}r^{\alpha -1}\mathrm{e}^{-ur^{\alpha }}du.
\end{align*}%
\noindent Inverting the Laplace transform yields 
\begin{equation*}
D_{t}^{\alpha }\hat{p}_{s}^{\ast
}(t,v)=(s+s^{2}+...+s^{k}-k)\int\limits_{0}^{\infty }\lambda (u+v)\exp
\{\Lambda (v,u+v)(s+s^{2}+...+s^{k}-k)\}h_{\alpha }(t,u)du=
\end{equation*}%
\begin{equation*}
=\int\limits_{0}^{\infty }\lambda (u+v)(s+s^{2}+...+s^{k}-k)\hat{p}%
_{s}(u,v)h_{\alpha }(t,u)du,
\end{equation*}%
where the mgf 
\begin{equation*}
\hat{p}_{s}(u,v)=\sum\limits\limits_{m}s^{m}p_{m}(u,v).
\end{equation*}
\noindent Finally, by inverting the mgf $(s+s^{2}+...+s^{k}-k)\hat{p}%
_{s}(u,v),$ we obtain: 
\begin{equation*}
D_{t}^{\alpha }p_{m}^{\ast }(u,v)=\int\limits_{0}^{\infty }\lambda (u+v)
\left [-kp_{m}(u,v)+\sum\limits_{j=1}^{m\wedge k}p_{m-j}(u,v) \right]
h_{\alpha }(t,u)du,
\end{equation*}%
since the mgf of 
\begin{equation*}
-kp_{m}(u,v)+\sum\limits_{j=1}^{m\wedge k}p_{m-j}(u,v)
\end{equation*}%
is equal to%
\begin{equation*}
\sum\limits_{m}s^{m}[-kp_{m}(u,v)+\sum\limits_{j=1}^{m\wedge
k}p_{m-j}(u,v)]= (s+s^{2}+...+s^{k}-k)\hat{p}_{s}(u,v).\qed
\end{equation*}
\noindent \textbf{Covariance structure}\newline
\noindent One can show that for NPPk \quad $\mathbb{E[}N^{n}(t)]=\frac{k(k+1)%
}{2}\Lambda (t),$ and its covariance function is 
\begin{equation*}
\mathrm{Cov}[N^{n}(t),N^{n}(s)]=\frac{k(k+1)(2k+1)}{6}\Lambda (\min (s,t)).
\end{equation*}%
\noindent Then the mean and covariance function of FNPPk are given by%
\begin{equation*}
\mathbb{E}[N_{\alpha }^{\ast }(t)]=\frac{k(k+1)}{2}\mathbb{E[}\Lambda
(Y_{\alpha }(t)\mathbb{]}
\end{equation*}%
\begin{equation*}
\mathrm{Cov}[N_{\alpha }^{\ast }(t),N_{\alpha }^{\ast }(s)]=\frac{%
k(k+1)(2k+1)}{6}\mathbb{E}\left[\Lambda (Y_{\alpha }(\min (s,t)) +\left ( \frac{%
k(k+1)}{2} \right)^{2}\mathrm{Cov}[\Lambda (Y_{\alpha }(t)),\Lambda
(N_{\alpha }(s))\right].
\end{equation*}


\section{P\'olya-Aeppli process of order $k$}

The P\'{o}lya-Aeppli process of order $k$ was defined and studied in the
context of ruin problems in \citep{Chuko2015} and later by \citep{kost2019}.
Related pure fractional birth processes were studied in \citep{Orsinger}.

\begin{definition}
The process $N_{PAk}(t)$ is said to be the P\'{o}lya-Aeppli process of order 
$k$ (PAk) if 
\begin{equation*}
N_{PAk}(t)=X_{1}+\dots +X_{N_{1}(t)},
\end{equation*}%
where (i) the random variables $X_{i}$ are i.i.d with the truncated
geometric distribution of parameter $\varrho \in \lbrack 0,1),$ given by (%
\ref{TG}); (ii) $N=\{N(t);t\geq 0\}$ is a homogeneous Poisson process (HPP)
with intensity $\lambda >0,$ independent of $\{X_{i}\}_{i=1}^{\infty }.$
\end{definition}

\noindent \noindent The following Kolmogorov forward equations are valid for
the marginal distribuions $p_{m}(t)=\mathbb{P}[N_{PAk}(t)=m]:$ 
\begin{align}
& \frac{d}{dt}p_{0}(t)=-\lambda p_{0}(t)  \label{svertka} \\
& \frac{d}{dt}p_{m}(t)=-\lambda p_{m}(t)+\lambda \frac{1-\varrho }{1-\varrho
^{k}}\sum\limits_{j=1}^{m\wedge k}\varrho ^{j-1}p_{m-j}(t),  \notag
\end{align}%
where $p_{m}(0)=\delta _{m,0.}$

The marginal distributions of the PAk process are given by 
\begin{equation}
p_{m}(t):=\mathbb{P}[N_{PAk}(t)=m]=q_{m}(\lambda t),m=0,1,2,..,  \label{pamd}
\end{equation}%
where $q_{m}$ are given by (\ref{PAPMF}).

More explicit expressions for $p_{m}(t)$ can be found in \citep{Min}.
\noindent The expectation ad variance are as follows: 
\begin{align}
& \mathbb{E}[N_{PAk}(t)]=\lambda t\frac{1+\varrho +\dots +\varrho
^{k-1}-k\varrho ^{k}}{1-\rho ^{k}},  \label{pak_moment} \\
& \mathrm{Var}[N_{PAk}(t)]=\frac{\lambda t}{1-\varrho ^{k}}\left[ 1+3\varrho
+5\varrho ^{2}+\dots +(2k-1)\varrho ^{k-1}-k^{2}\varrho ^{k}\right] .  \notag
\end{align}%
Note that PAk process is a compound Poisson process with the pgf 
\begin{equation*}
G_{N_{PAk}(t)}(u)=\mathbb{E}\left[ u^{N_{PAk}(t)}\right] =\mathbb{P}%
[N_{PAk}(t)=m]=e^{-\lambda t(1-G_{X}(u))},
\end{equation*}%
where $G_{X}(u)=\mathbb{E}\left[ u^{X}\right] $ is given by (\ref{pgf}%
).\noindent

\subsection{\textit{Non-homogeneous P\'olya-Aeppli process of order $k$}}

\noindent We now consider a non-homogeneous version by introducing a
deterministic time dependent intensity function $\lambda (t)$ as above, and $%
\Lambda (s,s+t)=\Lambda (s+t)-\Lambda (s), \:\Lambda (t)=\int_{0}^{t}\lambda
(u)du.$

\begin{definition}
(Non-homogeneous P\'olya-Aeppli process of order $k $). We define a
non-homogeneous P\'olya-Aeppli process of order $k $ with cumulative rate
function $\Lambda(t) $ and parameter $\varrho$ as 
\begin{align}
N^n_{PAk} (t) = X_1+ \dots + X_{N_1^n(t)},
\end{align}
where (i) $\{N_{1}^{n}(t); \: t\geq 0\}$ is a non-homogeneous Poisson
process (NPP) with cumulative rate function $\Lambda (t)$; (ii) $X_{i}$ are
i.i.d. r.v's following the truncated geometric distribution with parameter $%
\varrho,$ given by (\ref{pgf}); (iii) $\{ N_{1}^n(t); \: t \ge 0\}$ is
independent from $X_i, i=1,2,\dots$
\end{definition}

\noindent Note, that the random variable $%
N_{PAk}^{n}(t+s)-N_{PAk}^{n}(s),:s,:t\geq 0$ has the P\'{o}lya -Aeppli
distribution of order $k$ with parameters $\Lambda (s,t),\varrho ,$ that is%
\begin{equation}
f_{m}^{n}(t,u)=\mathbb{P}[N_{PAk}^{n}(t+u)-N_{PAk}^{n}(u)=m]=q_{m}(\Lambda
(u,u+t)),m=0,1,2,...,  \label{DPA}
\end{equation}%
where $q_{m}$ are given by (\ref{PAPMF}).

Then the marginal distributions of the process $N_{PAk}^{n}(t)$ are $\mathbb{%
P}[N_{PAk}^{n}(t)=m]=f_{m}^{n}(t,0).$ \noindent An alternative definition
can be given in terms of transition probabilities.

\begin{definition}
The counting process $N^n_{PAk} (t) $ is said to be a non-homogeneous
P\'olya -Aeppli process of order $k $ with the rate function $\lambda(t) $
and parameter $\varrho \in [0,1)$ if\newline
(1) $N^n_{PAk} (0)=0;$ (2) $N^n_{PAk} (t)$ has independent increments; (3)
for all $t \ge 0 $ 
\begin{equation}  \label{eqn: npak }
\mathbb{P}[N_{PAk}^n (t+h)=n \:|\:N_{PAk}^n (t)=m]= 
\begin{cases}
1-\lambda(t+h) h + o(h), \quad n=m \\ 
\frac{1-\varrho}{1-\varrho^k} \varrho^{i-1} \lambda(t+h) h + o(h), \quad
n=m+i, i= 1,2,\dots, k \\ 
\end{cases}%
\end{equation}
\end{definition}

\noindent It is easy to verify that the previous two definitions are
equivalent. \newline
\noindent \textbf{Marginal distributions of the process. }\newline
\noindent\ \noindent 
The following theorem holds.

\begin{theorem}
The functions $f_m(t,u), \: m=0,1,2, \dots $ satisfy the differential
equation: 
\begin{equation}  \label{eqn:  NPA}
\frac{d}{dt} f_m(t,u) = -\lambda(t+u) f_m(t,u) +\lambda(t+u) \frac{1-\varrho%
}{1-\varrho^k} \sum\limits_{j=1}^{m\wedge k }\varrho^{j-1} f_{x-j}(t,u).
\end{equation}
\end{theorem}

\noindent \textit{\ Proof.} We first consider the case $m=0.$ By fixing $u$
and taking a small $h$ we can write 
\begin{align*}
& f_{0}(t+h,u)=P[I(t+h)=0]=\mathbb{P}[N_{PAk}^{n}(t+u+h)-N_{PAk}^{n}(u)=0]=
\\
& \mathbb{P}[\mbox{no events in }:(u,u+t]:\cap \mbox{ no events in}%
:(u+t,u+t+h]]= \\
& \mathbb{P}[\mbox{ no  events in }:(u,u+t+h]]\mathbb{P}[\mbox{no events in}%
(u+t,u+t+h]= \\
& f_{0}(t+h)[1-\lambda (t+u)h+o(h)] \\
&
\end{align*}%
Thus 
\begin{equation*}
\frac{f_{0}(t+h,u)-f_{0}(t,u)}{h}=-\lambda (t+u)f_{0}(t,u)+\frac{o(h)}{h}.
\end{equation*}%
Letting $h\rightarrow 0$ yields 
\begin{equation*}
\frac{d}{dt}f_{0}(t,u)=-\lambda (t+u)f_{0}(t,u).
\end{equation*}%
For $m\geq 1$ we have 
\begin{align*}
& f_{m}(t+h,u)=\mathbb{P}[\{\mbox{ m  events in }:(u,u+t+h]\cap \{%
\mbox{no
events in}(u+t,u+t+h]\} \\
& :\cup \{m-1\mbox{  events in }:(u,u+t]\}\cap \mbox{1 event in}%
:(u+t,u+t+h]\}\cup \dots \\
& \cup \{\mbox{ 0 events in }:(u,u+t]\cap \mbox{m events in}(u+t,u+t+h]\}]=
\\
& f_{m}(t+h,u) \left[ 1-\lambda (t+u)h+o(h)]+f_{m-1}(t+h,u)[\frac{1-\varrho 
}{1-\varrho ^{k}}\lambda (t+u)h\varrho ^{1-1}+o(h) \right]+... \\
& +f_{0}(t+h,u)\left[\frac{1-\varrho }{1-\varrho ^{k}}\lambda (t+u)h\varrho
^{m\wedge k-1}+o(h)\right]= \\
& =\lambda (t+u)f_{m}(t+h,u)+\lambda (t+u)\frac{1-\varrho }{1-\varrho ^{k}}%
\sum\limits_{j=1}^{m\wedge k}\varrho ^{j-1}f_{m-j}(t,u). \\
\end{align*}%
%
%
%
%
%
%
\noindent Letting $h\rightarrow 0$ yields 
\begin{equation*}
\frac{d}{dt}f_{m}(t,u)=-\lambda (t+u)f_{m}(t,u)+\lambda (t+u)\frac{1-\varrho 
}{1-\varrho ^{k}}\sum\limits_{j=1}^{m\wedge k}\varrho ^{j-1}f_{m-j}(t,u),
\end{equation*}%
which was the statement of the theorem. \qed\newline
Note that in case $k\rightarrow \infty $, the P\'{o}lya-Aeppli process $%
N_{PAk}^{n}(t)$ coincides with the non-homogeneous P\'{o}lya-Aeppli process
defined in \citep{chuko2018}, but for a fixed $k$ the the P\'{o}lya-Aeppli
process $N_{PAk}^{n}(t)$ is new. 

\subsection{\textit{Fractional P\'olya-Aeppli process of order $k$}}

To the best of our knowledge, fractional versions of PAk processes have not
been considered yet. We define a fractional P\'olya-Aeppli process of order $%
k $ as a P\'olya-Aeppli process of order $k $ time-changed by the process $%
\{Y_{\alpha} (t); \:t \ge 0\},$ such that 
\begin{align}
N_{\alpha} ^{h}(t) =N_{PAk}(Y_{\alpha}(t)), \quad 0<\alpha<1,
\end{align}
where (i) $N_{1}=\{N_{1}(t); \: t \geq 0\}$ is the homogeneous Poisson
process with intensity $\lambda;$ (ii) $N_{PAk}(t)=X_{1}+\dots
+X_{N_{1}(t)}; $ (iii) $\{Y_{\alpha }(t); \: t\geq 0\},0<\alpha <1$ is the
inverse $\alpha $-stable subordinator, defined in (\ref{inverse}) and
independent of $N_{1}(t).$\newline
\noindent \textbf{Marginal distributions}\newline
\noindent We shall now obtain governing equations for the marginal
distributions of the fractional PAk process 
\begin{equation*}
p_{x}^{\alpha }(t)=\mathbb{P}[N_{PAk}(Y_{\alpha
}(t))=m]=\int\limits_{0}^{\infty }p_{m}(u)h_{\alpha }(t,u)du,\quad
m=0,1,\dots ,
\end{equation*}%
where $p_{m}(u)$ is given by (\ref{pamd}).

\begin{theorem}
The probabilities $p_{x}^{\alpha }(t),\quad x=0,1,\dots $ satisfy the
fractional differential-difference equations: 
\begin{align}
& D_{t}^{\alpha }p_{0}^{\alpha }(t)=-\lambda p_{0}^{\alpha }(t) \\
& D_{t}^{\alpha }p_{x}^{\alpha }(t)=-\lambda p_{x}^{\alpha }(t)+\lambda 
\frac{1-\varrho }{1-\varrho ^{k}}\sum\limits_{j=1}^{x\wedge k}\varrho
^{j-1}p_{x-j}^{\alpha }(t),
\end{align}%
where $D_{t}^{\alpha }f(t)$ is the fractional Caputo derivative of the
function $f$ given by (\ref{caputo}).
\end{theorem}

\noindent \textit{\ Proof.} \noindent We first consider the case $m\geq 1.$
By taking the fractional Caputo derivative of the both sides (\ref{svertka})
and using the property (\ref{governinvsub}), we get 
\begin{align*}
D_{t}^{\alpha }p_{m}^{\alpha }(t)& =-\int\limits_{0}^{\infty }p_{m}(u)\frac{%
\partial }{\partial u}h_{\alpha }(t,u)du= \\
& =\int\limits_{0}^{\infty } \left[ -\lambda p_{m}(u)+\lambda \frac{%
1-\varrho }{1-\varrho ^{k}}\sum\limits_{j=1}^{x\wedge k}\varrho
^{j-1}p_{m-j}(t) \right] h_{\alpha }(t,u)du-p_{m}(u)h_{\alpha
}(t,u)|_{0}^{\infty }= \\
& =-\lambda p_{m}^{\alpha }(t)+\lambda \frac{1-\varrho }{1-\varrho ^{k}}%
\sum\limits_{j=1}^{x\wedge k}\varrho ^{j-1}p_{m-j}^{\alpha }(t).
\end{align*}%
%
%
%
%
%
%
%
%
%
%
%
%
%
%
For $m=0$ we have 
\begin{align*}
D_{t}^{\alpha }p_{0}^{\alpha }(t)& =-\int\limits_{0}^{\infty }p_{0}(u)\frac{%
\partial }{\partial u}h_{\alpha }(t,u)du= \\
& =\int\limits_{0}^{\infty }[-\lambda p_{0}(u)]h_{\alpha }(t,u)du=-\lambda
p_{0}^{\alpha }(t).\qed
\end{align*}

\subsection{\textit{Correlation structure and long-range dependence property}%
}

In this sub-section we shall obtain several important characteristics of the
fractional P\'{o}lya-Aeppli process of order $k$ such as its expectation,
variance and covariance. After that, we are able to study the correlation
structure of the process. For the fractional P\'{o}lya-Aeppli process of
order $k$, $N_{\alpha }^{h}(t)=N_{PAk}(Y_{\alpha }(t)),$ we can use the
property of the conditional expectation to write (see \citep{Leo1}) 
\begin{align*}
& \mathbb{E}[N_{\alpha }^{h}(t)]=\mathbb{E}[\mathbb{E}[N_{\alpha }^{h}(t)\:|\:Y_{\alpha }(t)]\:|\:Y_{\alpha }(t)]=\int\limits_{0}^{\infty }%
\mathbb{E}[N_{PAk}(u)]h_{\alpha }(t,u)du= \\
& =\lambda \mathbb{E}[N_{PAk}(1)]\frac{t^{\alpha }}{\Gamma (\alpha +1)},
\end{align*}%
\begin{equation*}
\mathrm{Var}[N_{\alpha }^{h}(t)]=\frac{t^{\alpha }\mathrm{Var}%
[N_{PAk}(1)]}{\Gamma (\alpha +1)}+\frac{t^{2\alpha }(\mathbb{E}%
[N_{PAk}(1)])^{2}}{\alpha }\left( \frac{1}{\Gamma (2\alpha )}-\frac{1}{%
\alpha \Gamma (\alpha )^{2}}\right) .
\end{equation*}%
The covariance function can be calculated via the formula: 
\begin{equation*}
\mathrm{Cov}[N_{\alpha }^{h}(t),N_{\alpha }^{h}(s)]=\mathrm{Var}%
[N_{PAk }(1)]\frac{\min(t,s)^{\alpha }}{\Gamma (1+\alpha )}+(\mathbb{E}%
[N_{PAk }(1)])^{2}\mathrm{Cov}[Y_{\alpha }(t),Y_{\alpha }(s)],
\end{equation*}%
where the covariance of the process $Y_{\alpha }(t)$ is given by equation (%
\ref{covinv}).

\begin{theorem}
The process $N_{\alpha}^{h}(t) $ has the LRD property.
\end{theorem}

\noindent \textit{\ Proof.} Using the results from \citep{Leo1} similarly to
the previous section, we get 
\begin{align*}  \label{corr}
&\mathrm{Corr}[N^h_{\alpha}(t), N^h_{\alpha}(s)] \sim t^{-\alpha} C(\alpha, s)
\qquad t \to \infty,
\end{align*}
where $C(\alpha, s)= \left(\frac{1}{\Gamma(2\alpha)} - \frac{1}{
\alpha(\Gamma(\alpha))^2}\right)^{-1} \left[ \frac{\alpha \mathrm{Var}[N_{PAk}(1)] 
}{\Gamma(1+\alpha) (\mathbb{E}[N_{PAk}(1)] )^2} +\frac{ \alpha s^{\alpha}}{%
\Gamma(1+2\alpha)} \right]$, and $\mathbb{E}[N_{PAk}(1)]$,  $\mathrm{Var}[N_{PAk}(1)] $
are given by (\ref{pak_moment}). 
Thus the correlation function of FPAk process decays at rate $t^{-\alpha},
\: \alpha\in(0,1)$ and satisfies the LRD property. \qed

\subsection{\textit{\ Non-homogeneous fractional PAk process }}

As we did before, we can now define a non-homogeneous fractional P\'{o}%
lya-Aeppli process of order $k$ as 
\begin{equation*}
N_{\alpha }^{n}(t)=N_{PAk}(\Lambda (Y_{\alpha }(t))),\quad t\geq 0,\quad
0<\alpha <1,
\end{equation*}%
where all the symbols have the usual meaning defined above. We assume that
the inverse subordinator $Y_{\alpha }$ is independent of the process $%
N_{PAk}.$ In this sub-section, we shall derive governing equations for the
probabilities 
\begin{equation*}
p_{m}^{\ast \ast }(t,v)=\mathbb{P}[N_{PAk}(\Lambda (Y_{\alpha
}(t)+v))-N_{PAk}(\Lambda (v))=m].
\end{equation*}%
%
%
%
%
%
%
%
%
%
%
%
%
%
%

\begin{theorem}
The marginal distributions $p_{x}^{\ast \ast }(t,v)$ satisfy the following
fractional differential-difference integral equations 
\begin{align}
& D_{t}^{\alpha }p_{0}^{\ast \ast }(u,v)=-\int\limits_{0}^{\infty }\lambda
(u+v)f_{0}^{n}(u,v)h_{\alpha }(t,u)du  \label{govern_fnpak} \\
& D_{t}^{\alpha }p_{m}^{\ast \ast }(u,v)=\int\limits_{0}^{\infty }\lambda
(u+v)[-f_{m}^{n}(u,v))+\frac{1-\varrho }{1-\varrho ^{k}}\sum\limits_{j=1}^{m%
\wedge k}\rho ^{j-1}f_{m-j}^{n}(u,v)]h_{\alpha }(t,u)du\quad m=1,2,\dots
\end{align}%
with the initial condition $p_{m}^{\ast \ast }(0,v)=\delta _{m,0},$ where $%
f_{m}^{n}(u,v)$ is given by (\ref{DPA}).
\end{theorem}

\noindent \textit{\ Proof.} Using (\ref{LTIS}), the mgf of $f_{m}^{n}(u,v)$
can be written in the form: 
\begin{equation*}
\hat{f}_{s}^{n}(u,v)=\mathbb{E} \left[ s^{N^{n}(v+u)-N^{n}(v)} \right]=\exp
\left \{\Lambda (v,u+v)\frac{1-\rho }{1-\rho ^{k}}\sum\limits_{j=1}^{k}\rho
^{j-1}(s^{j}-1) \right \},
\end{equation*}%
while the Laplace transform with respect to $t$ of $h_{\alpha }(t,u)$ is
given by (\ref{LTIS}). Taking both the mgf and the Laplace transform in (\ref%
{Main}) as above, we have 
\begin{multline}
\bar{p}_{s}^{\ast \ast }(u,v)=r^{\alpha -1}\int\limits_{0}^{\infty }\hat{f}%
_{s}^{n}(u,v)\tilde{h}_{\alpha }(r,u)du= \\
\int\limits_{0}^{\infty } \left[ \exp \left \{\Lambda (v,u+v)\frac{1-\rho }{%
1-\rho ^{k}}\sum\limits_{j=1}^{k}\rho ^{j-1}(s^{j}-1) \right \} \right
]e^{-ur^{\alpha }}du.  \label{11}
\end{multline}
Note that for $U(u)=\exp \left \{\Lambda (v,u+v)\frac{1-\rho }{1-\rho ^{k}}%
\sum\limits_{j=1}^{k}\rho ^{j-1}(s^{j}-1) \right \}$ one can take derivative
in $u$ as follows: 
\begin{equation}
\frac{d}{du}U(u)=\frac{1-\rho }{1-\rho ^{k}}\sum \limits_{j=1}^{k}\rho
^{j-1}(s^{j}-1)[\lambda (v,u+v)]\exp \left \{\Lambda (v,u+v)\frac{1-\rho }{%
1-\rho ^{k}}\sum\limits_{j=1}^{k}\rho ^{j-1}(s^{j}-1) \right\}.  \label{12}
\end{equation}
\noindent Thus, integrating (\ref{11}) by parts with 
\begin{equation*}
U=\exp \left \{\Lambda (v,u+v)\frac{1-\rho }{1-\rho ^{k}}\sum%
\limits_{j=1}^{k}\rho ^{j-1}(s^{j}-1) \right\}, \quad V=-\frac{1}{r^{\alpha }%
}e^{-ur^{\alpha }},
\end{equation*}%
we get%
\begin{equation*}
\bar{p}_{s}^{\ast \ast }(u,v)=\frac{1}{r^{\alpha }} \left[ r^{\alpha -1}+%
\frac{1-\rho }{1-\rho ^{k}}[(s-1)+\rho (s^{2}-1)+...+\rho ^{k-1}(s^{k}-1) %
\right]\times
\end{equation*}%
\begin{equation}
\times \int\limits_{0}^{\infty }\lambda (v,u+v)\exp \left \{\Lambda (v,u+v)%
\frac{1-\rho }{1-\rho ^{k}}\sum\limits_{j=1}^{k}\rho ^{j-1}(s^{j}-1) \right
\}\mathrm{e}^{-ur^{\alpha }}r^{\alpha -1}du  \label{13}
\end{equation}
\noindent where $p_{s}^{\ast \ast }(0^{+},v)=1, $ since \: $Y_{\alpha }(0)=0$
a.s. Hence, by (\ref{13}) 
\begin{align*}
& r^{\alpha }\bar{p}_{s}^{\ast \ast }(r,v)-r^{\alpha -1}\bar{p}_{s}^{\ast
\ast }(0,v)=\mathcal{L}_{r}\{D_{t}^{\alpha }\bar{p}_{s}^{\ast \ast
}(r,v)\}(r)= \\
& =\frac{1-\rho }{1-\rho ^{k}}[(s-1)+\rho (s^{2}-1)+ \dots +\rho
^{k-1}(s^{k}-1)] \times \\
&\times \int\limits_{0}^{\infty }\lambda (u+v)\exp \left \{\Lambda (v,u+v)%
\frac{1-\rho }{1-\rho ^{k}}\sum\limits_{j=1}^{k}\rho ^{j-1}(s^{j}-1) \right
\} r^{\alpha -1}\mathrm{e}^{-ur^{\alpha }}du.
\end{align*}%
Inverting the Laplace transform yields 
\begin{align*}
& D_{t}^{\alpha }\hat{p}_{s}^{\ast \ast }(t,v)=\frac{1-\rho }{1-\rho ^{k}}%
[(s-1)+\rho (s^{2}-1)+...+\rho ^{k-1}(s^{k}-1)] \times \\
&\times \int\limits_{0}^{\infty }\lambda (u+v)\exp \left \{\Lambda (v,u+v)%
\frac{1-\rho }{1-\rho ^{k}}\sum\limits_{j=1}^{k}\rho ^{j-1}(s^{j}-1) \right
\}h_{\alpha }(t,u)du= \\
&=\int\limits_{0}^{\infty }\lambda (u+v) \left [\frac{1-\rho }{1-\rho ^{k}}%
[(s-1)+\rho (s^{2}-1)+\dots +\rho ^{k-1}(s^{k}-1) \right] \hat{f}%
_{s}^{n}(u,v)h_{\alpha }(t,u)du,
\end{align*}%
where the mgf is 
\begin{equation*}
\hat{f}_{s}^{n}(u,v)=\sum \limits_{m}s^{m}f_{m}^{n}(u,v).
\end{equation*}
Finally, by inverting the mgf 
\begin{equation*}
\left \lbrack \frac{1-\rho }{1-\rho ^{k}}[(s-1)+\rho (s^{2}-1)+...+\rho
^{k-1}(s^{k}-1) \right] \hat{f}_{s}^{n}(u,v),
\end{equation*}%
we obtain: 
\begin{equation*}
D_{t}^{\alpha }p_{m}^{\ast }(u,v)=\int\limits_{0}^{\infty }\lambda (u+v) 
\left[ -f_{m}^{n}(u,v)+\frac{1-\rho }{1-\rho ^{k}}\sum\limits_{j=1}^{m\wedge
k}\rho ^{j-1}f_{m-j}^{n}(u,v) \right] h_{\alpha }(t,u)du.\qed
\end{equation*}

\section{Discussion on simulations, applications and extensions}

The counting processes of order $k$ that we have discussed in this paper
have this general form 
\begin{equation}  \label{generalform}
N(t) = \sum_{i=1}^{\mathcal{N}(t)} X_i,
\end{equation}
where $\{X_i \}_{i=1}^\infty$ is a sequence of i.i.d. integer random
variables assuming values in $1, \dots, k$ and $\mathcal{N}(t)$ is a
counting process independent from the sequence. One further assumes that $%
N(0)=0$. Thanks to the independence among the increments and between them
and the counting process, one can express the moments of $N(t)$ as a
function of the moments of $X_1$ and $\mathcal{N}(t)$. In particular, one
can use Wald's equation for the expected value as all the assumptions of
Wald's theorem are satisfied \citep{wald44,wald45} 
\begin{equation}  \label{wald}
\mathbb{E}[N(t)] = \mathbb{E}[\mathcal{N}(t)] \mathbb{E}[X_1].
\end{equation}
Similarly, for the variance of $N(t)$, one gets 
\begin{equation}  \label{compoundvariance}
\mathrm{Var}[N(t)] = \mathbb{E}[\mathcal{N}(t)] \mathrm{Var}[X_1] + \mathrm{%
Var}[\mathcal{N}(t)] (\mathbb{E}[X_1])^2.
\end{equation}
Equation \eqref{generalform} suggests a simple and straightforward algorithm
for the simulation of $N(t)$. Given a value of $t$, one can first extract a
value $\mathcal{N}(t)=n$ and then sum $n$ values independently drawn from
the distribution of $X_1$. A simple algorithm in R that performs this task
is given in the appendix when $\mathcal{N}(t)$ is the fractional Poisson
process of renewal type used above and discussed by Mainardi et al. %
\citep{Mainardi1} and when $X_1$ is uniformly distributed in $1, \dots, k$.

\noindent This paper introduces a non-homogeneous fractional Poisson process
of order $k $ and a non-homogeneous fractional P\'olya-Aeppli process of
order $k. $ These processes lead to numerous generalizations and
applications. One of the natural extensions of this research is to
application in insurance where claims arrive in clusters in a non-homogenous
way. A natural feature of these processes is their long-range dependence
structure. These processes are therefore potentially useful for modelling
high-frequency financial data as well. Back to insurance, the generalized
classical ruin problems with FPPk and FPAk processes will be considered in
our next paper in the spirit of the results in \citep{Biard1}, \citep{Mahe1}%
, and \citep{kost2019}. A further extension would be considering a
time-change related to other non-local operators, rather than the Caputo
derivative, as investigated in e.g. \citep{toaldo}. This is also the subject
of ongoing research.

\noindent Kreer \cite{kreer} proves dynamical scaling under suitable
hypotheses for the non-homogeneous fractional Poisson process. The author
goes on to conjecture a mechanism according to which, the non-homogeneous
fractional Poisson process can be applied to phase transitions in which
nucleation is described by a counting process. It will be interesting to see
whether this potential application to the physics of phase transitions of
non-homogeneous fractional processes can lead to falsifiable experimental
predictions.

\section*{Acknowledgements}
The authors thank Mr. Mostafizar Khandakar and Dr. Kuldeep Kumar Kataria for their thorough reading of the manuscript that highlighted some typos and a wrong formula. Enrico Scalas acknowledges partial support from  the {\em Dr Perry James (Jim) Browne Research Centre} at the Department of Mathematics, University of Sussex.

\section*{Appendix}

The program listed below simulates the random variable $N(t)$ when the
process is the fractional Poisson process of order $k$. In order to generate
the time change, it is useful to generate Mittag-Leffler distributed i.i.d.
random variables $J_i$. A simple algorithm to do so is described in %
\citep{fulger} based on previous work on Linnik distributions %
\citep{kozubowski}. It gives the value of the random variable $J_1$ as 
\begin{equation*}
j_1 = - \gamma \log(u) \left( \frac{\sin(\alpha \pi)}{\tan(\alpha \pi v)} -
\cos(\alpha \pi)\right)^{1/\alpha}
\end{equation*}
where $\alpha \in (0,1)$ is the parameter of the Mittag-Leffler
distribution, $(u,v) \in (0,1)^2$ are independent uniform random variables
and $\gamma$ is a scale parameter. The code can be easily modified to
include any distribution on the integers $1, \ldots, k$ and any counting
process $\mathcal{N}(t)$.\newline

\texttt{\noindent \# Fractional Poisson process of order k }

\texttt{\noindent \# Parameters of the distribution }

\texttt{\noindent k <- 3 }

\texttt{\noindent lambda <- 1 }

\texttt{\noindent alpha <- 0.95 }

\texttt{\noindent \# Time for the process }

\texttt{\noindent t <- 10 }

\texttt{\noindent \# We compute Iter values of N(t) }

\texttt{\noindent Iter <- 10000 }

\texttt{\noindent Nt <- c() }

\texttt{\noindent for (i in 1:Iter)\{ }

\texttt{\noindent \# We generate a Fractional Poisson random variable }

\texttt{\noindent u <- runif(100000) }

\texttt{\noindent v <- runif(100000) }

\texttt{\noindent tau <-
-lambda*log(u)*(sin(alpha*pi)/tan(alpha*pi*v)-cos(alpha*pi))$^\wedge$%
(1/alpha) }

\texttt{\noindent time <- cumsum(tau) }

\texttt{\noindent n <- length(which(time<=t)) }

\texttt{\noindent \# We sum uniform random variables in 1, ..., k from 1 to
n }

\texttt{\noindent N <- 0 }

\texttt{\noindent if (n==0) \{Nt=c(Nt,0)\} }

\texttt{\noindent for (j in 1:n) \{ }

\texttt{\noindent X <- sample(k,1) }

\texttt{\noindent N <- N+X }

\texttt{\noindent \} }

\texttt{\noindent Nt <- c(Nt,N) }

\texttt{\noindent \} }


\end{document}